\documentclass[11pt,noinfoline]{article}
\RequirePackage[OT1]{fontenc}

\newcommand{\R}{\mathbb{R}}

\newcommand{\ave}{\mbox{ave}}

\newcommand{\eps}       {\varepsilon}
\newcommand{\PP} {{  \rm I\hskip-0.22em P}}

\newcommand{\EE} {{\rm I\hskip-0.48em E}}

\usepackage{graphicx}
\usepackage{multirow}
\usepackage{amsfonts}
\usepackage{natbib}
\usepackage{amsmath}

\newtheorem{Theo}{Theorem}
\newtheorem{Def}{Definition}
\newtheorem{Coro}{Corollary}
\newtheorem{Lemma}{Lemma}
\newtheorem{Proposition}{Proposition}
\newtheorem{Corr}{Corrolary}

\begin{document}
\bibliographystyle{plain}


\title{Estimating high-dimensional directed acyclic graphs
with the PC-algorithm}



\author{Markus Kalisch and Peter B\"uhlmann\footnote{Both authors are
    affiliated with the Seminar f\"ur Statistik, ETH Z\"urich, Switzerland.}}


\maketitle

\begin{abstract}
We consider the PC-algorithm (\cite{sgs}) for estimating the skeleton
of a very 
high-dimensional acyclic  
directed graph (DAG) with corresponding Gaussian distribution. The
PC-algorithm is computationally 
feasible for sparse problems with many nodes, i.e. variables, and it has the
attractive property to automatically achieve high computational
efficiency as a function of sparseness of the true underlying DAG.    
We prove consistency of the algorithm for very high-dimensional, sparse
DAGs where 
the number of nodes is allowed to quickly grow with sample size $n$, as fast as
$O(n^a)$ for any $0<a<\infty$. The sparseness assumption
is rather minimal requiring only that the neighborhoods in the
DAG are of lower order than sample size $n$. 
We empirically demonstrate the PC-algorithm for 
simulated data and argue that the algorithm is rather insensitive to
the choice of its single tuning parameter.    

\end{abstract}


\section{Introduction}
Graphical models  are a popular probabilistic
tool to analyze and visualize 
conditional independence relationships between random variables 
(see \cite{edwards}, \cite{lauritzen}). 
Major 
building blocks of the models are nodes, which represent random variables
and edges, which encode conditional dependence relations of the
enclosing vertices. The structure of conditional independence among
the random variables can be explored using the Markov properties.

Of particular current interest are directed acyclic graphs (DAGs),
containing directed rather than undirected edges, which restrict in a sense
the conditional dependence relations. These graphs can be interpreted by
applying the 
directed Markov property. When ignoring the directions of a DAG, we get the
skeleton of a DAG. In general, it is different from the conditional
independence graph (CIG), see section \ref{subsec.defi}. 
Thus, estimation methods for directed graphs cannot be easily borrowed from
approaches for undirected CIGs. 

Estimation of a DAG from data is difficult and computationally non-trivial
due to the enormous size of the space of DAGs: the number of possible DAGs
is super-exponential in the number of nodes. Nevertheless, there are
quite successful search-and-score methods for problems where the number of
nodes is small or moderate. For example, the search space may be restricted to
trees as in MWST (Maximum Weight Spanning Trees; see \cite{chow} and
\cite{heckerman}), or a greedy 
search is employed. The greedy DAG search can be
improved by exploiting probabilistic equivalence relations, and the search
space can be reduced from individual DAGs to equivalence classes, as
proposed in GES (Greedy 
Equivalent Search, see \cite{chickering}). Although this method seems
quite promising when having few or a moderate number of nodes only, it is
limited by the fact that the space of 
equivalence classes is conjectured to grow super-exponentially in the nodes
as well (see \cite{countEqClasses}). Bayesian approaches for DAGs,
which are computationally very intensive, include
\cite{spiegelhalter} and \cite{heckerman}.

An interesting alternative to greedy or structurally restricted approaches
is the PC-algorithm from \cite{sgs}. It starts from a
complete, undirected graph 
and deletes recursively edges based on conditional independence
decisions. This yields an undirected graph which can then be partially
directed and further extended to DAGs. For 
the skeleton of a DAG, i.e. the undirected version of a DAG, the
PC-algorithm runs in the worst case in exponential time (as a function of the
number of nodes), but if the true underlying DAG is sparse, which is often
a reasonable assumption, this reduces to a polynomial runtime. 

We focus in this paper on estimating DAGs in the high-dimensional context
when having many nodes, i.e. the number of nodes $p$ may be much larger than
sample size $n$. We prove that the PC-algorithm
consistently estimates the skeleton of an underlying sparse DAG, as sample
size $n \to
\infty$, even if $p = p_n = O(n^a)\ (0 < a < \infty)$ is allowed to grow
very quickly as a function of $n$. Our implementation of the PC-algorithm
allows to 
estimate the skeleton of a sparse DAG even if $p$ is in the hundreds or
thousands. For the high-dimensional setting with $p > n$, sparsity of the
underlying DAG is crucial for statistical consistency and computational
feasibility. The PC-algorithm seems to be the only method for high-dimensional
settings which is computationally feasible and, due to the new results in
this paper, provably correct in an asymptotic sense.  



We argue empirically that the PC-algorithm is rather insensitive to the
choice of its single tuning parameter, a significance level for testing,
and we compare 
the PC-algorithm with other methods, at least for low- or
mid-dimensional problems.

\section{The skeleton of a DAG}

\subsection{Definitions and preliminaries}\label{subsec.defi}
A graph $G = (V,E)$ consists of a set of nodes or vertices
$V=\{1,\ldots,p\}$ and a 
set of edges $E \subseteq V \times V$, i.e. the edge set is a subset of
ordered pairs of distinct nodes. In our setting, the set of nodes
corresponds to the components of a random vector $\mathbf{X} \in \R^p$. An edge
$(i,j) \in E$ is called directed if $(i,j) \in E$ but $(j,i) \notin E$: we
then use the notation $i \to j$. An acyclic directed graph (DAG) is a graph
$G$ where all edges are directed and not containing any cycle. 

If there is a directed edge $i \to j$, node $i$ is said to be a parent of
node $j$. The set of parents of node $j$ is denoted by $pa(j)$. The set of
neighbors of a node $j$, denoted by $ne(j)$, are all nodes $i$ with a
directed edge $i \to j$ or $j \to i$. Equivalently, $ne(j)$ is often
referred to as the adjacency set $adj(G,j)$ of a node $j$ in the graph
$G$. 
The skeleton of a DAG $G$ is
the undirected graph obtained from $G$ by substituting undirected edges
for directed edges.  

A probability distribution $P$ on $\R^p$ is said to be faithful with respect to
a graph $G$ if conditional independencies of the distribution can be
inferred from d-separation in the graph $G$ and vice-versa. More precisely:
consider a random vector $\mathbf{X} \sim P$. Faithfulness of $P$ with
respect to $G$ means: for every set $\mathbf{s} \subseteq V$, 
\begin{eqnarray*}
& &\mathbf{X}^{(i)}\ \mbox{and}\ \mathbf{X}^{(j)}\ \mbox{are conditionally
    independent given}\ \{\mathbf{X}^{(r)};\ r \in \mathbf{s}\}\\
&\Leftrightarrow & \mbox{node}\ i\ \mbox{and node}\ j\ \mbox{are
    d-separated by the set}\ \mathbf{s}.
\end{eqnarray*}
The notion of d-separation can be defined via moral graphs; details are
described in \cite[Prop. 3.25]{lauritzen}. We
    remark here that 
faithfulness is ruling out some classes of probability distributions. An
example of a non-faithful distribution is given in
\cite[Chapter 3.5.2]{sgs}. On the other hand, non-faithful
distributions form a 
Lebesgue 
null-set in the space of distributions associated with a DAG $G$, see
\cite[Th. 3.2]{sgs}. 

It is well known that for a probability distribution $P$ which is generated
from a DAG $G$, there is a whole equivalence class of DAGs with
corresponding distribution $P$ (see \cite[Section 2.2 ]{chickering}),
and we 
can only 
identify an equivalence class of DAGs, even when having infinitely many
observations. But the skeletons of DAGs from the same equivalence class
are the same, and thus, inferring a skeleton from data is an easier and
better identifiable task than aiming for directed graphs. 
We point out that
in general, the skeleton of a DAG $G$ with corresponding distribution $P$ is 
different from the conditional independence graph corresponding to the
distribution $P$. 
In particular, if $P$ is faithful with respect to a DAG $G$,
\begin{eqnarray}\label{skelchar}
& &\mbox{there is an edge between nodes $i$ and $j$ in the skeleton of DAG
  $G$}\nonumber\\ 
&\Leftrightarrow &\mbox{for all $\mathbf{s} \subseteq V \setminus \{i,j\}$},\
  \mathbf{X}^{(i)}\ \mbox{and}\ \mathbf{X}^{(j)}\ \mbox{are  
  conditionally independent}\nonumber\\
& & \mbox{\hfill given}\ \{\mathbf{X}^{(r)};\ r \in \mathbf{s}\},
\end{eqnarray} 
(\cite[Th. 3.4]{sgs}). This
implies the following: if $P$ is faithful with respect to a DAG $G$, the
skeleton of the DAG $G$ is a subset (or equal) to the conditional
independence graph (CIG) corresponding to $P$. The reason is that an edge
in a CIG requires only conditional dependence given the set $V \setminus
\{i,j\}$. 
We conclude that if the true underlying probability mechanisms are generated  
from a DAG, it
is more appropriate to use the undirected skeleton as a target than 
the undirected conditional independence graph.   

\subsection{The PC-algorithm for the skeleton}

A naive strategy would be to check conditional independencies given all
subsets $\mathbf{s} \subseteq V \setminus \{i,j\}$ (see formula
(\ref{skelchar})), 
i.e. all partial correlations in 
the case of multivariate normal distributions. This would become
computationally infeasible and statistically ill-posed for $p$ larger than
sample size.  A much better approach is to use the PC-algorithm which is
able to exploit sparseness of the  
graph. More precisely, we apply the part of
the PC-algorithm that identifies the undirected edges of the DAG. 

\subsubsection{Population Version}\label{subsec.pop}
In the population version of the PC-algorithm, we assume that perfect knowledge
about all necessary conditional independence relations is available. 
\begin{center}
\textbf{The PC$_{pop}$($m$)-algorithm}
\end{center}
\begin{enumerate}
\item Form the complete undirected graph $\tilde{C}$ on the vertex set V.
\item
Set $\ell = -1; \quad C=\tilde{C}$
\item[a)]repeat

Increase $\ell$ by one. 
\begin{enumerate}
\item[b)]repeat

Select an ordered pair of nodes $i$,$j$ that are adjacent in $C$ such that
$|adj(C,i) \setminus \{j\}| \geq \ell$ and $\mathbf{k} \subseteq adj(C,i)
\setminus \{j\}$ with $|\mathbf{k}|=\ell$. If $i$ and $j$ are conditionally
independent given $\mathbf{k}$, delete edge
$i,j$. Denote this new graph by $C$.
\item[b)]until all ordered pairs of adjacent variables $i$ and $j$ such
  that $|adj(C,i) 
\setminus \{j\}| \geq \ell$ and $\mathbf{k} \subseteq adj(C,i) \setminus
  \{j\}$ with 
$|\mathbf{k}|=\ell$ have been tested for conditional independence
\end{enumerate}
\item[a)]  until $\ell = m$ or\\
for each ordered pair of adjacent nodes $i$,$j$: $|adj(C,i) \setminus
\{j\}| \le \ell$.
\end{enumerate}
This is the description of the population PC$_{pop}$($m$)-algorithm which
is stopped at a pre-specified level 
$m$; the index $\ell$ may not even reach $m$ if the second statement for
termination of 2a) applies. There is no need to tune the parameter $m$ when
using the reached stopping level,
\begin{eqnarray}\label{reach}
m_{reach} =\ \max\{\mbox{stopping level}\ m;\ \mbox{index}\ \ell = m\}.
\end{eqnarray}
The value of $m_{reach}$ depends on the underlying distribution. 
\begin{Def}(Population version)
The  PC$_{pop}$-algorithm (\cite{sgs}) is
defined as the PC$_{pop}$($m_{reach}$)-algorithm. 
\end{Def}




A proof that this algorithm produces the correct skeleton 
can be easily deduced from Theorem 5.1 in \cite{sgs}. We summarize the
result as follows.  
\begin{Proposition} \label{theo:1}
Consider a DAG $G$ and assume that the distribution $P$ is faithful to
$G$. Denote the maximal number of neighbors by $q= \max_{1 \le j \le p}
|ne(j)|$. 
Then, the PC$_{pop}$-algorithm constructs the true skeleton of the DAG.   
Moreover, for the reached stopping level: $m_{reach} \in \{q-1,q\}$.  
\end{Proposition}

A proof is given in section \ref{sec.proofs}.

\subsubsection{Sample version for the skeleton}\label{subsec.sample}
For finite samples, we need to estimate conditional independencies. 
We limit ourselves to the Gaussian case, where all nodes correspond to
random variables with a multivariate normal distribution. Furthermore,
we assume faithful models, i.e. the
conditional independence relations can be read of the graph and vice
versa; see section \ref{subsec.defi}. 

In the Gaussian case, conditional
independencies can be inferred from partial correlations. 

\begin{Proposition} \label{lemma:1}
Assume that the distribution $P$  of the random vector $\mathbf{X}$ is
multivariate normal.
For $i \neq j \in \{1,\ldots ,p\},\ \mathbf{k} \subseteq \{1,\ldots ,p\}
  \setminus 
\{i,j\}$, denote by $\rho_{i,j|\mathbf{k}}$ the 
partial correlation between $\mathbf{X}^{(i)}$ and $\mathbf{X}^{(j)}$ given
$\{\mathbf{X}^{(r)};\ r \in \mathbf{k}\}$. Then, 
$\rho_{i,j|\mathbf{k}}=0$ if and only if $\mathbf{X}^{(i)}$ and $\mathbf{X}^{(j)}$
are conditionally independent given $\{\mathbf{X}^{(r)};\ r \in \mathbf{k}\}$. 
\end{Proposition}
Proof: The claim is an elementary property of the multivariate normal
distribution, cf. \cite[Prop. 5.2.]{lauritzen}. \hfill$\Box$ 

\medskip
We can thus estimate partial correlations to obtain estimates of  conditional
independencies. The sample partial correlation
$\hat{\rho}_{i,j|\mathbf{k}}$ can be 
calculated via regression or recursively by using the following identity:
for some $h \in \mathbf{k}$,
\begin{eqnarray*}
\rho_{i,j|\mathbf{k}} = \frac{\rho_{i,j|\mathbf{k} \setminus h} -
    \rho_{i,h| \mathbf{k} \setminus h} 
    \rho_{j,h|\mathbf{k} \setminus h}}{\sqrt{(1 - \rho_{i,h|\mathbf{k}
    \setminus h}^2)(1 - \rho_{j,h|\mathbf{k} \setminus h}^2)}}.
\end{eqnarray*}
For testing whether a partial correlation is zero or not, we apply Fisher's
z-transform
\begin{eqnarray}\label{ztrans}
Z(i,j|\mathbf{k}) = \frac{1}{2} \log \left (\frac{1 +
    \hat{\rho}_{i,j|\mathbf{k}}}{1 - \hat{\rho}_{i,j|\mathbf{k}}} \right).
\end{eqnarray}
Classical decision theory yields then the following rule when using the
significance level $\alpha$. Reject the
null-hypothesis $H_0(i,j|\mathbf{k}):\ \rho_{i,j|\mathbf{k}} = 0$ against
the two-sided 
alternative $H_A(i,j|\mathbf{k}):\ \rho_{i,j|\mathbf{k}} \neq 0$ if
$\sqrt{n-|\mathbf{k}|-3} |Z(i,j|\mathbf{k})| > \Phi^{-1}(1 
- \alpha/2)$, where $\Phi(\cdot)$ denotes the cdf of ${\cal N}(0,1)
$.   

The sample version of the PC-algorithm is almost identical to the
population version in section \ref{subsec.pop}, except from step 2b).

\begin{center}
\textbf{The PC($m$)-algorithm}
\end{center}
\begin{enumerate}
\item[]
Run the PC$_{pop}$($m$)-algorithm as described in section \ref{subsec.pop}
but replace in 
2b) the statement about conditional independence of $i,j$ given $\mathbf{k}$ by
$\sqrt{n-|\mathbf{k}|-3} |Z(i,j|\mathbf{k})| \le \Phi^{-1}(1 - \alpha/2)$,
see (\ref{ztrans}).   
\end{enumerate}

\medskip\noindent
The algorithm yields a data-dependent value $\hat{m}_{reach,n}$
which is the maximal stopping level that is reached, i.e. the sample version of
(\ref{reach}).
\begin{Def}(Sample version) The PC-algorithm is defined as the PC($\hat{m}_{reach,n}$)-algorithm. 
\end{Def}
As we will see in Theorem \ref{theo:2},
the stopping level $\hat{m}_{reach,n}$ provides a reasonable value for the
stopping level $m$. The only tuning parameter of the PC-algorithm is $\alpha$,
i.e. the significance level for testing partial correlations. The algorithm
seems to be rather insensitive to the choice of $\alpha$, see section
\ref{sec.numer}. 

As we
will see below in section \ref{sec.cons}, the algorithm is asymptotically
consistent even if 
$p$ is much larger than $n$ but the DAG is sparse.


\section{Consistency for high-dimensional skeletons}\label{sec.cons}  
We will show that the PC-algorithm from section \ref{subsec.sample} is
asymptotically consistent for the skeleton of a DAG, even if $p$ is much
larger than $n$ but the DAG is sparse. We assume that the data are realizations
of i.i.d. random vectors $\mathbf{X}_1,\ldots ,\mathbf{X}_n$ with
$\mathbf{X}_i \in \R^p$ from a DAG $G$ with corresponding distribution
$P$. To capture high-dimensional behavior, we will allow to let the
dimension grow as a function of sample size: thus, $p = p_n$ and also the
DAG $G = G_n$ and the distribution $P = P_n$.  
Our assumptions are as follows. 
\begin{enumerate}
\item[(A1)] The distribution $P_n$ is multivariate Gaussian and faithful to
  the DAG $G_n$ for all $n$. 
\item[(A2)] The dimension $p_n = O(n^a)$ for some $0 \le a < \infty$. 
\item[(A3)] The maximal number of neighbors in the DAG $G_n$ is denoted by\\
  $q_n = \max_{1 \le j \le p_n} |ne(j)|$, with $q_n = O(n^{1 - b})$ for
  some $0 < b \le 1$.  
\item[(A4)] The partial correlations between
    $\mathbf{X}^{(i)}$ and $\mathbf{X}^{(j)}$ given $\{\mathbf{X}^{(r)}; r
    \in \mathbf{k}\}$ for some set $\mathbf{k} \subseteq \{1,\ldots ,p_n\}
    \setminus \{i,j\}$ are denoted by $\rho_{n;i,j|\mathbf{k}}$. Their
    absolute values are
    bounded from below and above: 
\begin{eqnarray*}
& &\inf\{|\rho_{i,j|\mathbf{k}}|;\ i,j,\mathbf{k}\ \mbox{with}\
\rho_{i,j|\mathbf{k}} \neq 0\} \ge c_n,\ c_n^{-1} = O(n^{d}),\\ 
& &\hspace*{55mm}\ \mbox{for some}\ 0 < d < b/2,\\
& &\sup_{n;i,j,\mathbf{k}} |\rho_{i,j|\mathbf{k}}| \le M < 1,
\end{eqnarray*}
where $0 < b \le 1$ is as in (A3).
\end{enumerate}
Assumption (A1) is an often used assumption in graphical modeling, although
it does restrict the class of possible probability distributions (see also
third paragraph of section \ref{subsec.defi}); (A2) allows
for an arbitrary polynomial growth of dimension as a function 
of sample size, i.e. high-dimensionality; (A3) is a sparseness assumption and
(A4) is a regularity  
condition. Assumptions (A3) and (A4) are rather minimal:
note that 
with $b = 1$ in (A3), e.g. fixed $q_n =q < \infty$, $m_n = m < \infty$, the
partial 
correlations can decay as $n^{-1/2 + \eps}$ for any $0 < \eps \le
1/2$. Our assumptions are 
simpler and seem to be weaker, although not directly comparable, than in
\cite{mebu06} who 
analyze the Lasso for estimating high-dimensional undirected conditional
independence graphs (where the growth in dimensionality is as in (A2)). If
the dimension $p$ is fixed (with fixed DAG $G$ and fixed distribution $P$),
(A2), (A3) and (A4) hold and (A1) remains as the only condition. 

\begin{Theo}\label{theo:1a}
Assume (A1), (A2), (A3) with $0< b \le 1$ and (A4) with $0<d<b/2$. Denote
by $\hat{G}_{skel,n}(\alpha_n,m_n)$ the estimate 
from the PC($m_n$)-algorithm in section \ref{subsec.sample} and by
$G_{skel,n}$ the true skeleton from the DAG $G_n$. Moreover, denote by
$m_{reach,n}$ the value described in (\ref{reach}). Then, for $m_n \ge
m_{reach,n},\  m_n = O(n^{1-b})\ (n \to \infty)$, there exists $\alpha_n \to 0\ (n \to
\infty)$ such 
that  
\begin{eqnarray*}
& &\PP[\hat{G}_{skel,n}(\alpha_n,m_n) = G_{skel,n}] \\
&=& 1 - O(\exp(-Cn^{1 - 2d})) \to 1\ (n \to
\infty)\ \mbox{for some}\ 0<C < \infty.
\end{eqnarray*}
\end{Theo}

A proof is given in section \ref{sec.proofs}. The lower
bound of the range for $m_n$ is $m_{reach,n}$ is either equal to
$q_n-1$ or $q_n$, see Proposition \ref{theo:1}, i.e. it depends on the
unknown sparseness $q_n$ in 
(A3). A non-constructive choice for the value of the significance level is
$\alpha_n = 2(1 - \Phi(n^{1/2} c_n/2))$ which depends on the 
unknown lower bound of partial correlations in (A4).

\medskip\noindent
\textbf{Remark 1.} For the case with fixed dimension $p$ (with fixed DAG
$G$ and fixed distribution $P$) , Theorem \ref{theo:1a}  becomes: for any
choice of $m_n \ge 
p-2,\ m_n = o(n)\ (n \to \infty)$ and using $\alpha_n = 2(1 - \Phi(D (n
\log(n)^{-1})^{1/2}))$ for any $0<D<\infty$,
\begin{eqnarray*}
& &\PP[\hat{G}_{skel,n}(\alpha_n,m_n) = G_{skel}] \\
&=& 1 - O(\exp(-Cn \log(n)^{-1})) \to 1\ (n \to
\infty)\ \mbox{for some}\ 0<C < \infty.
\end{eqnarray*}

\medskip\noindent
\textbf{Remark 2.} Denote by $u_n$ the minimal stopping level $m$ such that
the population PC-algorithm PC$_{pop}$($m$) yields the true skeleton of
the underlying DAG $G$. It is known that $u_n \le \max_{1 \le j \le p_n}
|pa(j)|$, i.e. the maximal number of parents; this can be deduced from
Theorem 5.1 in \cite{sgs}. Moreover, Theorem \ref{theo:1a} also holds for
$m_n \ge u_n,\ m_n = O(n^{1-b})$, and instead of (A3) it
would suffice to require the weaker condition that $u_n = O(n^{1-b})$. The
latter holds if  
the maximal number of parents satisfies $\max_{1 \le j \le p_n}|pa(j)| =
O(n^{1-b})$. The proof is as for Theorem \ref{theo:1a}.  

\medskip
Theorem \ref{theo:1a} leaves some flexibility for choosing $m_n$. The
PC-algorithm yields a data-dependent reached stopping level
$\hat{m}_{reach,n}$, i.e. the sample version of (\ref{reach}). 
\begin{Theo}\label{theo:2}
Assume (A1)-(A4). Then,
\begin{eqnarray*}
& &\PP[\hat{m}_{reach,n} = m_{reach,n}] = 1 - O(\exp(-Cn^{1-2d})) \to 1\
(n \to \infty)\\
& &\mbox{for some}\ 0<C < \infty,
\end{eqnarray*}
where $d>0$ is as in (A4). 
\end{Theo}
 
A proof is given in section \ref{sec.proofs}. 
Because there are faithful
distributions which require $m_n = m_{reach,n} \in \{q_n-1,q_n\}$ for
consistent estimation with the PC($m$)-algorithm,
Theorem \ref{theo:2} indicates that the PC-algorithm,
stopping at $\hat{m}_{reach,n}$, yields with high probability the smallest
$m = m_n$ which 
is universally consistent for all faithful distributions.
Therefore, there
is no need to select a tuning parameter $m = m_n$: the PC-algorithm yields
a good, data-dependent $\hat{m}_{reach,n}$.  

Theorems \ref{theo:1a} and \ref{theo:2} together yield the
consistency of the PC-algorithm, i.e. the PC($\hat{m}_{reach,n}$)-algorithm.  
\begin{Corr}
Assume (A1)-(A4). Denote by $\hat{G}_{skel,n}(\alpha_n)$ the estimate
from the PC-algorithm in section \ref{subsec.sample} and by
$G_{skel,n}$ the true skeleton from the DAG $G_n$. Then, there exists
$\alpha_n \to 0\ (n \to \infty)$ such that  
\begin{eqnarray*}
& &\PP[\hat{G}_{skel,n}(\alpha_n) = G_{skel,n}] \\
&=& 1 - O(\exp(-Cn^{1 - 2d})) \to 1\ (n \to
\infty)\ \mbox{for some}\ 0<C < \infty,
\end{eqnarray*}
where $d >0$ is as in (A4). 
\end{Corr}

Our theoretical framework allows for rather large values of $p$. The
computational complexity of the PC-algorithm is difficult to evaluate
exactly, but the worst case is bounded by  
\begin{eqnarray}\label{comp}
O(p^{\hat{m}_{reach,n}})\ \mbox{which is with high probability bounded by}\
O(p^{q_n}) 
\end{eqnarray}
as a function of dimensionality $p$. We note that the bound may be very
loose for many distributions. Thus, for the worst case where the complexity
bound is achieved, the algorithm is computationally feasible if $q_n$ is
small, say $q_n \le 3$, even if $p$ is large. For non-worst cases, however,
we can still do the computations for much larger values of $q_n$ and fairly
dense graphs, e.g. some nodes $j$ have neighborhoods of size up to $|ne(j)|
= 30$. 
 
In practice, we can check the value of $\hat{m}_{reach,n}$. As long as it is
of ``lower order'' than sample size $n$, the PC-algorithm
yields satisfactory results.

\section{Numerical examples}\label{sec.numer}
We analyze the PC-algorithm and other alternative methods for the skeleton
using various simulated data.
The numerical results have been obtained using the \texttt{R}-package
\texttt{pcalg} (\cite{ka05}) and the \texttt{Bayes Net Toolbox} of
Kevin Murphy. 

\subsection{Simulating data}\label{subsec.simmod}
In this section, we analyze the PC-algorithm for the skeleton using
simulated data.

In order to simulate data, we first construct an adjacency matrix $A$ as
follows: 
\begin{enumerate}
\item Fix an ordering of the variables.
\item Fill the adjacency matrix $A$ with zeros.
\item Replace every matrix entry in the lower triangle (below the diagonal) by
  independent realizations of Bernoulli($s$) random variables with success
  probability $s$ where $0<s<1$. We will call $s$
  the sparseness of the model. 
\item Replace each entry with a 1 in the adjacency matrix by independent
  realizations of a Uniform($[0.1,1]$) random variable.
\end{enumerate}
This then yields a matrix $A$ whose entries are zero or in the range
$[0.1,1]$. The corresponding DAG draws a directed edge from node $i$ to
node $j$ if $i < j$ and $A_{ji} \neq 0$. 
The DAGs (and skeletons thereof) that are created in this way have the
following property: $\EE{[N_i]} = s(p-1)$, where $N_i$ is the number of
neighbors of a node $i$.

Thus, a low sparseness parameter $s$ implies few neighbors and vice-versa. 
The matrix $A$ will be used to generate the data as follows. The
value of the random variable $X^{(1)}$, corresponding to the first node, is
given by  
\begin{eqnarray*}
& &\epsilon^{(1)} \sim N(0,1) \\
& &X^{(1)} = \epsilon^{(1)}
\end{eqnarray*}
and the values of the next random variables (corresponding to the next
nodes) can be computed recursively as
\begin{eqnarray*}
& &\epsilon^{(i)} \sim N(0,1) \\
& &X^{(i)} = \sum_{k=1}^{i-1} A_{ik} X^{(k)} + \epsilon^{(i)}\ (i=2,\ldots ,p),
\end{eqnarray*}
where all $\epsilon^{(1)},\ldots ,\epsilon^{(p)}$ are independent. 

\subsection{Comparison with alternative methods}
In this section, we will compare the PC-algorithm with two alternative
methods, Greedy Equivalent Search (GES, see \cite{chickering}) and Maximum
Weight Spanning Trees (MWST, see \cite{heckerman}) which both try to find
DAGs that maximize the BIC criterion. 

We found, that the BIC based methods find DAGs with high True
Positive Rate (TPR) but also rather high False
Positive Rate (FPR). If only a small amount of observations is available
(as is often the case in a very high-dimensional setting), we cannot hope to
recover the complete underlying model. Therefore, instead of
large TPR, we would rather prefer a subset of edges
with high reliability. 
A measure for high reliability is the True
\begin{table}[!htbp]
  \centering
\begin{tabular}{|r|r|r|r|}
\hline
Method & $\ave[TPR]$ & $\ave[FPR]$ & $\ave[TDR]$ \\
\hline
PC & 0.57 (0.06) & 0.02 (0.01) & 0.91 (0.05) \\
GES & 0.85 (0.05) & 0.13 (0.04) & 0.71 (0.07) \\
MWST & 0.66 (0.07) & 0.06 (0.01) & 0.78 (0.06) \\
\hline
\end{tabular}
\caption{$p=10$ nodes, sample size $n=50$, sparseness 
  $s = 0.1$, 50 replicates. Standard errors are given in parentheses.
   The PC-algorithm achieves a substantially higher True Discovery
       Rate than GES or MWST.}
\label{table:compareMethods}
\end{table}
Discovery Rate (TDR), which is the ratio of correctly found edges and the
total number of all edges found. 

As can be seen in table
\ref{table:compareMethods}, the PC-algorithm achieves in 
our simulations by far higher True Discovery Rates than GES or MWST: of all
found edges, 91\% were correct. Thus, although a smaller total of
edges was found, the estimated edges were correct more frequently. We think,
that this is a substantial advantage for real world applications.

\subsection{Different parameter settings}
As introduced in section \ref{subsec.sample}, the PC-algorithm has only one
tuning parameter $\alpha$. 
In this section, we analyze the dependence of the
algorithm on this parameter for different settings. 
\begin{table}[htbp]
  \centering
\begin{tabular}{|r|r|r|r|r|}
\hline
$\alpha$ & $\ave[TPR]$ & $\ave[FPR]$ & $\ave[TDR]$ & $\ave[\hat{m}_{reach}]$\\
\hline
0.001 & 0.065 (0.002) & 0.0057 (0.0005) & 0.80 (0.02) & 2.56 (0.07)\\
0.01 & 0.089 (0.003) & 0.0082 (0.0007) & 0.78 (0.02) & 2.92 (0.06)\\
0.05 & 0.116 (0.003) & 0.0133 (0.0009) & 0.75 (0.02) & 3.26 (0.06)\\
0.1 & 0.128 (0.003) & 0.0161 (0.0010) & 0.73 (0.02) & 3.46 (0.08)\\
0.3 & 0.151 (0.005) & 0.0238 (0.0011) & 0.68 (0.02) & 4.28 (0.08)\\
\hline
\end{tabular}
\caption{$p=30$, $n=20$, $s=0.1$, $50$ replicates; s.e. in parentheses.} 
\label{tab:n20s01}
\end{table}

\begin{table}[htbp]
  \centering
\begin{tabular}{|r|r|r|r|r|}
\hline
$\alpha$ & $\ave[TPR]$ & $\ave[FPR]$ & $\ave[TDR]$ & $\ave[\hat{m}_{reach}]$\\
\hline
0.001 & 0.069 (0.002) & 0.0056 (0.0005) & 0.80 (0.02) & 2.30 (0.07) \\
0.01 & 0.092 (0.002) & 0.0097 (0.0007) & 0.77 (0.02) & 2.92 (0.06) \\
0.05 & 0.116 (0.003) & 0.0141 (0.0008) & 0.73 (0.01) & 3.28 (0.07) \\
0.1 & 0.131 (0.003) & 0.0165 (0.0008) & 0.73 (0.01) & 3.50 (0.08) \\
0.3 & 0.159 (0.004) & 0.0233 (0.0010) & 0.70 (0.01) & 4.34 (0.07) \\
\hline
\end{tabular}
\caption{$p=30$, $n=20$, $s=0.4$, $50$ replicates; s.e. in parentheses.}
\label{tab:n20s04}
\end{table}

\begin{table}[htbp]
  \centering
\begin{tabular}{|r|r|r|r|r|}
\hline
$\alpha$ & $\ave[TPR]$ & $\ave[FPR]$ & $\ave[TDR]$ & $\ave[\hat{m}_{reach}]$ \\
\hline
0.001 & 0.153 (0.004) & 0.015 (0.001) & 0.77 (0.01) & 4.02 (0.07) \\
0.01 & 0.175 (0.005) & 0.017 (0.001) & 0.77 (0.01) & 4.38 (0.09) \\
0.05 & 0.193 (0.005) & 0.020 (0.001) & 0.76 (0.01) & 4.82 (0.08) \\
0.1 & 0.200 (0.005) & 0.021 (0.001) & 0.76 (0.01) & 5.00 (0.09) \\
0.3 & 0.221 (0.006) & 0.025 (0.001) & 0.74 (0.01) & 5.66 (0.09) \\
\hline
\end{tabular}
\caption{$p=30$, $n=100$, $s=0.1$, $50$ replicates; s.e. in parentheses.}
\label{tab:n100s01}
\end{table}

\begin{table}[htbp]
  \centering
\begin{tabular}{|r|r|r|r|r|}
\hline
$\alpha$ & $\ave[TPR]$ & $\ave[FPR]$ & $\ave[TDR]$ & $\ave[\hat{m}_{reach}]$\\
\hline
0.001 & 0.155 (0.004) & 0.015 (0.001) & 0.78 (0.01) & 4.12 (0.08) \\
0.01 & 0.174 (0.004) & 0.016 (0.001) & 0.78 (0.01) & 4.54 (0.08) \\
0.05 & 0.188 (0.005) & 0.020 (0.001) & 0.76 (0.01) & 4.78 (0.09) \\
0.1 & 0.196 (0.005) & 0.021 (0.001) & 0.76 (0.01) & 4.92 (0.09) \\
0.3 & 0.217 (0.006) & 0.028 (0.001) & 0.71 (0.01) & 5.58 (0.10) \\
\hline
\end{tabular}
\caption{$p=30$, $n=100$, $s=0.4$, $50$ replicates; s.e. in parentheses.}
\label{tab:n100s04}
\end{table}

\begin{table}[htbp]
  \centering
\begin{tabular}{|r|r|r|r|r|}
\hline
$\alpha$ & $\ave[TPR]$ & $\ave[FPR]$ & $\ave[TDR]$ & $\ave[\hat{m}_{reach}]$ \\
\hline
0.001 & 0.250 (0.007) & 0.033 (0.001) & 0.71 (0.01) & 6.5 (0.1) \\
0.01 & 0.258 (0.007) & 0.036 (0.001) & 0.70 (0.01) & 6.7 (0.1) \\
0.05 & 0.264 (0.007) & 0.038 (0.001) & 0.69 (0.01) & 7.0 (0.1) \\
0.1 & 0.268 (0.007) & 0.041 (0.001) & 0.68 (0.01) & 7.3 (0.1) \\
0.3 & 0.283 (0.007) & 0.047 (0.001) & 0.67 (0.01) & 7.6 (0.1) \\
\hline
\end{tabular}
\caption{$p=30$, $n=5000$, $s=0.1$, $50$ replicates; s.e. in parentheses.}
\label{tab:n5000s01}
\end{table}

\begin{table}[htbp]
  \centering
\begin{tabular}{|r|r|r|r|r|}
\hline
$\alpha$ & $\ave[TPR]$ & $\ave[FPR]$ & $\ave[TDR]$ & $\ave[\hat{m}_{reach}]$ \\
\hline
0.001 & 0.260 (0.007) & 0.031 (0.001) & 0.73 (0.01) & 6.40 (0.09) \\
0.01 & 0.268 (0.007) & 0.035 (0.001) & 0.72 (0.01) & 6.80 (0.09)\\
0.05 & 0.277 (0.006) & 0.036 (0.001) & 0.72 (0.01) & 7.04 (0.09)\\
0.1 & 0.281 (0.007) & 0.038 (0.001) & 0.71 (0.01) & 7.22 (0.10)\\
0.3 & 0.294 (0.006) & 0.045 (0.001) & 0.68 (0.01) & 7.70 (0.11)\\
\hline
\end{tabular}
\caption{$p=30$, $n=5000$, $s=0.4$, $50$ replicates; s.e. in parentheses.}
\label{tab:n5000s04}
\end{table}

Tables \ref{tab:n20s01} to \ref{tab:n5000s04} show the average over 50
replicates 
of TPR, FPR, TDR and $\hat{m}_{reach}$ for the DAG model in section
\ref{subsec.simmod} with $p = 30$ nodes and varying sample size $n$ and
sparseness $s$. 

In the wide range of $\alpha$s, no choice can be identified as being the
best or worst. Especially in the case of very few observations
we see that small $\alpha$ leads to the discovery of very few edges with
high reliability (high TDR), whereas higher values of $\alpha$ lead to the
discovery of more edges but with less reliability. Therefore, $\alpha$ can
be used for fine tuning in finding a good compromise between amount of
edges found and their reliability. 

Note, however, that especially for
larger sample sizes, the rates vary only little, sometimes only by a few
percent. Comparing this with the large change in $\alpha$ (over two orders
of magnitude), we feel that the PC-algorithm is rather insensitive to the
choice of its single tuning parameter.

\section{Conclusions}
The PC-algorithm is a powerful method for estimating the skeleton of a
potentially very high-dimensional DAG with corresponding Gaussian
distribution. Sparsity, in terms of the maximal size of the neighborhoods
of the true underlying DAG, is crucial for statistical 
consistency (assumption (A3) and Theorem 
\ref{theo:1a}) and for computational feasibility with at most a polynomial
complexity (see (\ref{comp})) as a function of dimensionality. We prove
consistency for high-dimensional frameworks under rather minimal assumption
on sparseness and decay of non-zero partial correlations.   



The PC-algorithm compares well with alternative
approaches like MWST and GES for low- or mid-dimensional problems. For
high-dimensional settings, MWST and GES (with the implementations we used)
become extremely slow while the 
PC-algorithm is still computationally feasible; e.g. a polynomial algorithm
for a sparse DAG, see (\ref{comp}). Software for the PC-algorithm will be
made available in \texttt{R}, package \texttt{pcalg} (\cite{ka05}). 

%
%
%

\section{Proofs}\label{sec.proofs}

\subsection{Proof of Proposition \ref{theo:1}}
Consider $\mathbf{X}$ with distribution $P$. Since $P$ is faithful
to the DAG $G$, conditional independence of $\mathbf{X}^{(i)}$ and
$\mathbf{X}^{(j)}$ given $\{\mathbf{X}^{(r)};\ r \in \mathbf{k}\}$
($\mathbf{k} \subseteq V \setminus\{i,j\}$) is equivalent to d-separation of
nodes $i$ and $j$ given 
the set $\mathbf{k}$ (see \cite[Th. 3.3]{sgs}). Thus, the population
PC$_{pop}$-algorithm as formulated in section \ref{subsec.pop} coincides with the 
one from \cite{sgs} which is using the concept of d-separation, and the
first claim about correctness of the skeleton follows from \cite[Th. 5.1.,
Ch. 13]{sgs}. 

The second claim about the value of $m_{reach}$ can be proved as
follows. First, due to the definition of the PC$_{pop}$($m$)-algorithm and
the fact that it constructs the correct skeleton,  $m_{reach} \le q$. We
now argue that $m_{reach} \ge q-1$. Suppose the contrary. Then, $m_{reach}
\le q-2$: we could then continue with a further iteration in the
algorithm since $m_{reach} + 1 \le q-1$ and there is at least one node $j$
with neighborhood-size $|ne(j)| = q$: that is, the reached stopping level
would be at least $q-1$ which is a contradiction to $m_{reach} \le
q-2$.\hfill$\Box$ 

\subsection{Proof of Theorem \ref{theo:1a}}

\subsubsection{Analysis of partial correlations}
We first establish uniform consistency of estimated partial
correlations. Denote by $\hat{\rho}_{i,j}$ and $\rho_{i,j}$ the sample and
population correlation between $\mathbf{X}^{(i)}$ and
$\mathbf{X}^{(j)}$. Likewise, $\hat{\rho}_{i,j|\mathbf{k}}$ and
$\rho_{i,j|\mathbf{k}}$ denote the sample and population partial
correlation between $\mathbf{X}^{(i)}$ and $\mathbf{X}^{(j)}$ given
$\{\mathbf{X}^{(r)}; r \in \mathbf{k}\}$, where $\mathbf{k} \subseteq
\{1,\ldots ,p_n\} \setminus \{i,j\}$. 

\medskip
Many partial correlations (and non-partial correlations) are tested for
being zero during the run of the PC($m_n$)-algorithm. For a fixed ordered
pair of nodes $i,j$, the conditioning sets are elements of  
\begin{eqnarray*}
K_{i,j}^{m_n} = \{\mathbf{k} \subseteq \{1,\ldots, p_n\} \setminus \{i,j\}:
|\mathbf{k}| \leq m_n\}
\end{eqnarray*}
whose cardinality is bounded by 
\begin{eqnarray}\label{add2}
|K_{i,j}^{m_n}| \le B p_n^{m_n}\ \mbox{for some $0< B < \infty$}.
\end{eqnarray}
\begin{Lemma} \label{lemma:2}
Assume (A1) (without requiring faithfulness) and $\sup_{n,i \neq j}
|\rho_{n;i,j}| \le M < 1$ (compare with (A4)). Then, for any $0< \gamma \le 2$,
\begin{displaymath}
\sup_{i,j,\mathbf{k} \in K_{i,j}^{m_n}}
\PP[|\hat{\rho}_{n;i,j}-\rho_{n;i,j}|>\gamma] 
\le C_1 (n-2) \exp \left((n-4) \log(\frac{4-\gamma^2}{4+\gamma^2})\right),
\end{displaymath}
for some constant $0< C_1 <\infty$ depending on $M$ only.
\end{Lemma}
Proof: We make substantial use of \cite{hotelling53}'s work. Denote
by $f_n(r,\rho)$ the probability density function of 
the sample correlation $\hat{\rho} = \hat{\rho}_{n+1;i,j}$ based on $n+1$
observations and by $\rho = \rho_{n+1;i,j}$ the population correlation. (It
is notationally easier to work with sample size $n+1$; and we just use the
abbreviated notations with $\hat{\rho}$ and $\rho$). 
For $0 < \gamma \le 2$, 
\begin{eqnarray*}
\PP[|\hat{\rho}-\rho|>\gamma] =\PP[\hat{\rho}<\rho -\gamma] +
\PP[\hat{\rho}>\rho+\gamma]. 
\end{eqnarray*}
It can be shown, that $f_n(r,\rho) = f_n(-r,-\rho)$, see
\cite[p.201]{hotelling53}. This symmetry implies,
\begin{eqnarray}\label{symmet}
\PP_{\rho}[\hat{\rho}<\rho -\gamma] = \PP_{\tilde{\rho}}[\hat{\rho} >
\tilde{\rho} + \gamma]\ \mbox{with}\ \tilde{\rho} = -\rho.
\end{eqnarray}
Thus, it suffices to show that $\PP[\hat{\rho}>\rho +\gamma] =
\PP_{\rho}[\hat{\rho}>\rho +\gamma] $ decays 
exponentially in $n$, uniformly for all $\rho$.  

It has been shown (\cite[p.201, formula (29)]{hotelling53}), that
for $-1< \rho < 1$,
\begin{eqnarray}\label{ADDNEW1}
\PP[\hat{\rho}>\rho + \gamma] \le \frac{(n-1) \Gamma(n)}{\sqrt{2 \pi}
    \Gamma{(n + 
    \frac{1}{2})}} M_0(\rho+\gamma) (1+\frac{2}{1-|\rho|})
\end{eqnarray}
with 
\begin{eqnarray}\label{ADDNEW2}
& &M_0(\rho + \gamma) =  \int_{\rho + \gamma}^1 (1-\rho^2)^{\frac{n}{2}}
(1-x^2)^{\frac{n-3}{2}} (1-\rho x)^{-n+\frac{1}{2}} d x \nonumber \\
&= & \int_{\rho+\gamma}^{1} (1-\rho^2)^{\frac{\tilde{n}+3}{2}}
(1-x^2)^\frac{\tilde{n}}{2} (1-\rho x)^{-\tilde{n}-\frac{5}{2}} d x\ \
(\mbox{using $\tilde{n} = n-3$})\nonumber \\
&\le & \frac{(1-\rho^2)^\frac{3}{2}}{(1-|\rho|)^{\frac{5}{2}}}
\int_{\rho+\gamma}^{1}(\frac{\sqrt{1-\rho^2} \sqrt{1-x^2}}{1-\rho
  x})^{\tilde{n}} dx \nonumber \\
&\le & \frac{(1-\rho^2)^\frac{3}{2}}{(1-|\rho|)^\frac{5}{2}} 2
\max_{\rho+\gamma \leq x \leq 1} (\frac{\sqrt{1-\rho^2}
  \sqrt{1-x^2}}{1-\rho x})^{\tilde{n}}.
\end{eqnarray}
We will show now that $g_{\rho}(x) = \frac{\sqrt{1-\rho^2}
  \sqrt{1-x^2}}{1-\rho x} < 1$ for all $\rho + \gamma \le x \le 1$ and $-1
  < \rho < 
1$ (in fact, $\rho \le 1 - \gamma$ due to the first restriction). Consider 
\begin{eqnarray}\label{ADDNEW3}
\sup_{-1<\rho<1;\rho+\gamma \le x \le 1} g_{\rho}(x) &=&
  \sup_{-1<\rho \le 1-\gamma} \frac{\sqrt{1-\rho^2} 
  \sqrt{1-(\rho+\gamma)^2}}{1-\rho(\rho+\gamma)} \nonumber \\
&= &\frac{\sqrt{1-\frac{\gamma^2}{4}}
  \sqrt{1-\frac{\gamma^2}{4}}}{1-(\frac{-\gamma}{2})(\frac{\gamma}{2})} =
  \frac{4-\gamma^2}{4+\gamma^2} < 1\ \mbox{for all $0 < \gamma \le 2$}.
\end{eqnarray}

Therefore, for $-1 < -M \leq \rho \leq M < 1$ (see assumption (A4)) and
using (\ref{ADDNEW1})-(\ref{ADDNEW3}) together with the fact that
$\frac{\Gamma(n)}{\Gamma(n+\frac{1}{2})} \leq const.$ with 
respect to $n$, we have
\begin{eqnarray*}
& &\PP[\hat{\rho}>\rho+\gamma]\\
&\le & \frac{(n-1) \Gamma(n)}{\sqrt{2 \pi}
  \Gamma(n+\frac{1}{2})}
  \frac{(1-\rho^2)^{\frac{3}{2}}}{(1-|\rho|)^{\frac{5}{2}}} 2
  (\frac{4-\gamma^2}{4+\gamma^2})^{\tilde{n}} (1+\frac{2}{1-|\rho|})\\
&\le & \frac{(n-1) \Gamma(n)}{\sqrt{2 \pi} 
  \Gamma(n+\frac{1}{2})} \frac{1}{(1-M)^{\frac{5}{2}}} 2
  (\frac{4-\gamma^2}{4+\gamma^2})^{\tilde{n}} (1+\frac{2}{1-M}) \leq{} \\
&\le &C_1 (n-1) (\frac{4-\gamma^2}{4+\gamma^2})^{\tilde{n}} = C_1 (n-1)
\exp((n-3) \log (\frac{4-\gamma^2}{4+\gamma^2})), \\ 
\end{eqnarray*}
where $0<C_1 <\infty$ depends on $M$ only, but not on $\rho$ or
$\gamma$. By invoking (\ref{symmet}), the proof is complete (note that the
proof assumed sample size $n+1$).\hfill$\Box$ 

\medskip
Lemma \ref{lemma:2} can be easily extended to partial correlations, as shown
by \cite{fisher24}, using projections for Gaussian distributions.

\begin{Lemma} \label{lemma:3} (Fisher, 1924)\\
Assume (A1) (without requiring faithfulness). 
If the cumulative distribution function of $\hat{\rho}_{n;i,j}$ is denoted by
$F(\cdot|n,\rho_{n;i,j})$, then the cdf 
of the sample partial correlation $\hat{\rho}_{n;i,j|\mathbf{k}}$ with
$|\mathbf{k}|=m < n-1$ is $F[\cdot|n-m,\rho_{n;i,j|\mathbf{k}}]$. That is, the
effective sample size is reduced by $m$. 
\end{Lemma}
A proof can be found in \cite{fisher24}; see also \cite{anderson84}.
\hfill$\Box$   

\medskip
Lemma \ref{lemma:2} and \ref{lemma:3} yield then the following. 
\begin{Coro}\label{coro:1a}
Assume (the first part of) (A1) and (the upper bound in) (A4). Then, for
any $\gamma >0$,
\begin{eqnarray*}
& &\sup_{i,j,\mathbf{k} \in
  K_{i,j}^{m_n}}\PP[|\hat{\rho}_{n;i,j|\mathbf{k}}-\rho_{n;i,j|\mathbf{k}}|>
  \gamma]\\  
&\le& C_1 (n-2 -m_n) \exp 
\left((n-4-m_n) \log(\frac{4-\gamma^2}{4+\gamma^2})\right),
\end{eqnarray*}
for some constant $0< C_1 < \infty$ depending on $M$ from (A4) only. 
\end{Coro}

The PC-algorithm is testing partial correlations after the
z-transform $g(\rho) = 0.5 \log((1+\rho)/(1- \rho))$. Denote by
$Z_{n;i,j|\mathbf{k}} = g(\hat{\rho}_{n;i,j|\mathbf{k}})$ and by
$z_{n;i,j|\mathbf{k}} = g(\rho_{n;i,j|\mathbf{k}})$. 
\begin{Lemma}\label{lemm:3a}
Assume the conditions from Corollary \ref{coro:1a}.
Then, for any $\gamma >0$,
\begin{eqnarray*}
& &\sup_{i,j,\mathbf{k} \in K_{i,j}^{m_n}} 
\PP[|Z_{n;i,j|\mathbf{k}}-z_{n;i,j|\mathbf{k}}|>\gamma] \\
&\le& O(n-m_n) \left(\exp((n-4 -m_n)
  \log(\frac{4-(\gamma/L)^2}{4+(\gamma/L)^2})) + \exp(-C_2(n-m_n))\right)  
\end{eqnarray*}
for some constant $0< C_2< \infty$ and $L = 1/(1 - (1+M)^2/4)$. 
\end{Lemma}
Proof: A Taylor expansion of the z-transform $g(\rho) = 0.5
\log((1+\rho)/(1- \rho))$ yields:
\begin{eqnarray}\label{taylor}
Z_{n;i,j|\mathbf{k}} - z_{n;i,j|\mathbf{k}} =
g'(\tilde{\rho}_{n;i,j|\mathbf{k}}) (\hat{\rho}_{n;i,j|\mathbf{k}} -
\rho_{n;i,j|\mathbf{k}}),
\end{eqnarray}
where $|\tilde{\rho}_{n;i,j|\mathbf{k}} - \rho_{n;i,j|\mathbf{k}}| \le
|\hat{\rho}_{n;i,j|\mathbf{k}} - \rho_{n;i,j|\mathbf{k}}|$. Moreover,
$g'(\rho) = 1/(1 - \rho^2)$. By applying Corollary \ref{coro:1a} with
$\gamma = \kappa = (1 - M)/2$ we 
have
\begin{eqnarray}\label{addnew2}
& &\sup_{i,j,\mathbf{k} \in K_{i,j}^{m_n}}
\PP[|\tilde{\rho}_{n;i,j|\mathbf{k}} - \rho_{n;i,j|\mathbf{k}}| \le
\kappa]\nonumber \\
&> & 1 - C_1 (n-2 -m_n) \exp (-C_2 (n - m_n)).  
\end{eqnarray}
Since 
\begin{eqnarray*}
& &g'(\tilde{\rho}_{n;i,j|\mathbf{k}}) = \frac{1}{1 -
\tilde{\rho}_{n;i,j|\mathbf{k}}^2} = \frac{1}{1 -
(\rho_{n;i,j|\mathbf{k}} + (\tilde{\rho}_{n;i,j|\mathbf{k}} -
\rho_{n;i,j|\mathbf{k}}))^2}\\
&\le& \frac{1}{1 - (M + \kappa)^2}\ \mbox{if}\
|\tilde{\rho}_{n;i,j|\mathbf{k}} - \rho_{n;i,j|\mathbf{k}}| \le \kappa,
\end{eqnarray*}
where we also invoke (the second part of) assumption (A4) for the last
inequality. Therefore, since $\kappa = (1 - M)/2$ yielding $1/(1 - (M
+ \kappa)^2) = L$, and using (\ref{addnew2}), we get 
\begin{eqnarray}\label{add4}
& &\sup_{i,j,\mathbf{k} \in
  K_{i,j}^{m_n}}\PP[|g'(\tilde{\rho}_{n;i,j|\mathbf{k}})| \le L] \nonumber \\ 
&\ge&  1 - C_1 (n-2 -m_n) \exp (-C_2 (n - m_n)).
\end{eqnarray}
Since $|g'(\rho)| \ge 1$ for all $\rho$, we obtain with (\ref{taylor}):

\begin{eqnarray}\label{add3}
& &\sup_{i,j,\mathbf{k} \in K_{i,j}^{m_n}}
\PP[|Z_{n;i,j|\mathbf{k}}-z_{n;i,j|\mathbf{k}}|>\gamma] \\ \nonumber
&\le & \sup_{i,j,\mathbf{k} \in K_{i,j}^{m_n}} 
\PP[|g'(\tilde{\rho}_{n;i,j|\mathbf{k}})| > L] 
+ \sup_{i,j,\mathbf{k} \in K_{i,j}^{m_n}}
\PP[|\hat{\rho}_{n;i,j|\mathbf{k}} - \rho_{n;i,j|\mathbf{k}}| > \gamma/L].  
\end{eqnarray}

Formula (\ref{add3}) follows from elementary probability calculations: for
two random 
variables $U,V$ with $|U| \ge 1$ ($|U|$ corresponding to 
$|g'(\tilde{\rho})|$ and $|V|$ to the difference $|\hat{\rho} - \rho|$),
\begin{eqnarray*}
\PP[|U V| > \gamma]& =&
\PP[|U V| > \gamma, |U| > L] + \PP[|U V| > \gamma, 1 \le |U| \le L] \\
&\le& \PP[|U| > L] + \PP[|V| > \gamma/L]. 
\end{eqnarray*}
The statement then follows from (\ref{add3}), (\ref{add4}) and Corollary
\ref{coro:1a}.\hfill$\Box$

\subsubsection{Analysis of the PC($m$)-algorithm}
The population version PC$_{pop}$($m_n$)-algorithm when stopped at
level $m_n = 
m_{reach,n}$ constructs the true skeleton according to Proposition
\ref{theo:1}. Moreover, the PC$_{pop}$($m$)-algorithm remains to be correct
when using $m \ge m_{reach,n}$.
An error
occurs in the sample PC-algorithm if there is a pair of nodes $i,j$ and a
conditioning set $\mathbf{k} \in K_{i,j}^{m_n}$ (although the algorithm is
typically only going through a random subset of $K_{i,j}^{m_n}$) where an
error event $E_{i,j|\mathbf{k}}$ occurs; $E_{i,j,\mathbf{k}}$ denotes that
``an error occurred when testing partial correlation for zero at nodes $i,j$
with conditioning set $\mathbf{k}$''. Thus,
\begin{eqnarray}\label{err-pc}
& &\PP[\mbox{an error occurs in the PC($m_n$)-algorithm}]\nonumber \\
&\le&  
P[\bigcup_{i,j,\mathbf{k} \in K_{ij}^{m_n}} E_{i,j|\mathbf{k}}]
\le O(p_n^{m_n+2}) \sup_{i,j,\mathbf{k}\in K_{ij}^{m_n}}
\PP[E_{i,j|\mathbf{k}}],
\end{eqnarray}
using that the cardinality of the set $|\{i,j,\mathbf{k}\in K_{ij}^{m_n}\}| =
O(p_n^{m_n+2})$, see also formula (\ref{add2}). Now
\begin{eqnarray}\label{ADD1}
E_{i,j|\mathbf{k}} = E_{i,j|\mathbf{k}}^I \cup E_{i,j|\mathbf{k}}^{II},
\end{eqnarray}
where 
\begin{eqnarray*}
& &\mbox{type I error}\ E_{i,j|\mathbf{k}}^I:\ \sqrt{n - |k| - 3}
   |Z_{i,j|\mathbf{k}}| > \Phi^{-1} (1- \alpha/2)\ \mbox{and}\
   z_{i,j|\mathbf{k}} =0,\\ 
& &\mbox{type II error}\ E_{i,j|\mathbf{k}}^{II}:\ \sqrt{n - |k| - 3}
   |Z_{i,j|\mathbf{k}}| \le \Phi^{-1}
   (1- \alpha/2)\ \mbox{and}\ z_{i,j|\mathbf{k}} \neq 0.
\end{eqnarray*}
Choose $\alpha = \alpha_n = 2(1 - \Phi(n^{1/2} c_n/2))$, where $c_n$ is from
(A4). Then,
\begin{eqnarray}\label{ADD2}
\sup_{i,j,\mathbf{k} \in K_{i,j}^{m_n}}  \PP[E_{i,j|\mathbf{k}}^I] & = &
\sup_{i,j,\mathbf{k} \in K_{i,j}^{m_n}}  \PP[|Z_{i,j|\mathbf{k}} -
z_{i,j|\mathbf{k}}| > (n/(n - |\mathbf{k}| - 3))^{1/2} c_n/2]\nonumber\\
&\le & O(n - m_n) \exp(-C_3 (n - m_n) c_n^2),
\end{eqnarray}
for some $0<C_3< \infty$ using Lemma \ref{lemm:3a} and the fact that
$\log(\frac{4-\delta^2}{4+\delta^2}) \sim -\delta^2/2$ as $\delta \to 0$.  
Furthermore, with the choice of $\alpha = \alpha_n$ above,
\begin{eqnarray*}
& &\sup_{i,j,\mathbf{k} \in K_{i,j}^{m_n}}  \PP[E_{i,j|\mathbf{k}}^{II}] =
\sup_{i,j,\mathbf{k} \in K_{i,j}^{m_n}} \PP[|Z_{i,j|\mathbf{k}}| \le \sqrt{n/(n
- |\mathbf{k}| - 3)} c_n/2] \\
&\le& \sup_{i,j,\mathbf{k} \in K_{i,j}^{m_n}} \PP[|Z_{i,j|\mathbf{k}} -
z_{i,j|\mathbf{k}}| > c_n(1 - \sqrt{n/(n - |\mathbf{k}| - 3)}/2)],
\end{eqnarray*}
because $\inf_{i,j;\mathbf{k} \in K_{i,j}^{m_n}} |z_{i,j|\mathbf{k}}| \ge
c_n$ since $|g(\rho)| \ge |\rho|$ for all $\rho$ and using assumption
(A4). By invoking Lemma \ref{lemm:3a} we then obtain:
\begin{eqnarray}\label{ADD3}
\sup_{i,j,\mathbf{k} \in K_{i,j}^{m_n}}  \PP[E_{i,j|\mathbf{k}}^{II}] \le
O(n - m_n) \exp(-C_4 (n - m_n) c_n^2)
\end{eqnarray}
for some $0<C_4 < \infty$. 
Now, by (\ref{err-pc})-(\ref{ADD3}) we get
\begin{eqnarray*}
& &\PP[\mbox{an error occurs in the PC($m_n$)-algorithm}]\\
&\le& O(p_n^{m_n+2} (n - m_m) \exp(-C_5 (n - m_n) c_n^2)) \\
&\le& O(n^{a (m_n+2)
  +1} \exp(-C_5 (n - m_n) n^{-2d})) \\
&= &O\left(\exp \left(a(m_n+2) \log(n) + \log(n) -
C_5 (n^{1 - 2d} - m_n n^{-2d})\right)\right) = o(1),
\end{eqnarray*}
because $n^{1 - 2d}$ dominates all other terms in the argument of the
$\exp$-function due to the assumption in (A4) that $d<b/2$. This completes
the proof.\hfill$\Box$  

\subsection{Proof of Theorem \ref{theo:2}}
 
Consider the population algorithm PC$_{pop}$($m$): the reached stopping
level satisfies $m_{reach} \in \{ q_n-1,q_n\}$, see Proposition
\ref{theo:1}. The sample PC($m_n$)-algorithm
with stopping level in the range of $m_{reach} \le m_n = O(n^{1-b})$,
coincides 
with the population version on a set $A$ having probability $P[A] = 1 -
O(\exp(-Cn^{1-2d}))$, see the last formula in the proof of Theorem
\ref{theo:1a}. Hence, on the set $A$, $\hat{m}_{reach,n} = m_{reach} \in
\{q_n-1,q_n\}$. The claim then follows from Theorem
\ref{theo:1a}.\hfill$\Box$     
 

\bibliography{biblio}

\end{document}